\documentstyle[12pt]{article}
\def\ba{\begin{eqnarray}}
\def\ea{\end{eqnarray}}

\def\lb{\label}
\def\be{\begin{equation}}
\def\ee{\end{equation}}

\begin{document}

\begin{center}
{\Large \bf Representations of A-type Hecke algebras}\footnote{The work
of the first author (A. I.) was supported by the grants
INTAS 03-51-3350 and RFBR 05-01-01086-a; the work of the second author
(O. O.)was supported by the ANR project GIMP No. ANR-05-BLAN-0029-01.}
\end{center}

\vspace{.5cm}

\begin{center}
\large{A.P.\,Isaev${}^{a}$ and O.\,Ogievetsky${}^{b}$ }
\end{center}

\begin{center}
${}^{a}$ Bogoliubov Laboratory of Theoretical Physics,
Joint Institute for Nuclear Research, \\
Dubna, Moscow region 141980, Russia \\
E-mail: isaevap@theor.jinr.ru
\\
\vspace{0.3cm}
${}^{b}$ Center of Theoretical Physics\footnote{Unit\'e Mixte de Recherche
(UMR 6207) du CNRS et des Universit\'es Aix--Marseille I,
Aix--Marseille II et du Sud Toulon -- Var; laboratoire affili\'e \`a la
FRUMAM (FR 2291)}, Luminy,
13288 Marseille, France \\
and P. N. Lebedev Physical Institute, Theoretical Department,
Leninsky pr. 53, 117924 Moscow, Russia \\
E-mail: oleg@cpt.univ-mrs.fr
\end{center}

\vspace{2cm}

{\bf Abstract.} We review some facts about the
representation theory of the Hecke algebra.
We adapt for the Hecke algebra case
the approach of \cite{OV} which was developed for the representation
theory of symmetric groups.
We justify an explicit construction
of the idempotents in the Hecke algebra in terms of Jucys-Murphy elements.
Ocneanu's traces for these idempotents
(which can be interpreted
as q-dimensions of corresponding irreducible representations
of quantum linear groups) are presented.

\section{Introduction}
\setcounter{equation}0

Main statements of
the representation theory of Hecke algebras are known
mostly due to the works by  V.Jones, I.V.Cherednik, G.Murphy,
R.Dipper and G.James, H.Wenzl, a.o.
(see, e.g., \cite{Jon1} -- \cite{W1}).
In this report
the approach of \cite{OV}, developed
for the representation theory of symmetric groups,
is generalized to the case of the $A$-type Hecke algebras.
Certain propositions below are given
without proofs due to
lack of space and, also, because
the corresponding statements for Hecke algebras
are proved like those for symmetric groups.

The importance of the theory of the A-type Hecke algebra $H_M$ is that
$H_M$ is the centralizer of the action of general linear quantum
groups $U_q(gl(N))$ in the tensor powers $V^{\otimes M}$ of the
vector representation $V$ of $U_q(gl(N))$.
We have shown recently \cite{IsOg} that an arbitrary representation of the
Hecke algebra $H_M$ defines
an integrable model on a chain with $M$ sites. This fact demonstrates the importance of the
representation theory of the Hecke algebra in the theory of integrable models also.

\section{$A$-Type Hecke algebras and Jucys - Murphy elements}
\setcounter{equation}0

A braid group ${\cal B}_{M+1}$ is generated by Artin elements $\sigma_i$ $(i=1, \dots M)$
subject to relations:
\be
\label{braidg}
\sigma_i \, \sigma_{i+1} \, \sigma_i =
\sigma_{i+1} \, \sigma_i \,  \sigma_{i+1} \; , \;\;\;
\sigma_{i} \,  \sigma_{j} = \sigma_{j} \,  \sigma_{i} \;\;
{\rm for} \;\; |i-j| > 1 \; .
\ee

An $A$-Type Hecke algebra $H_{M+1}(q)$ (see e.g. \cite{Jon1} and Refs. therein) is
a quotient of the group algebra of the braid group ${\cal B}_{M+1}$
by an additional relation
\be
\label{ahecke}
\sigma^2_i - 1 = (q -q^{-1}) \,  \sigma_i \; , \;\; (i = 1, \dots , M) \; .
\ee
Here $q \in {\bf C} \backslash \{0\}$ is a parameter.
The group algebra of ${\cal B}_{M+1}$ (\ref{braidg})
has an infinite dimension while its quotient $H_{M+1}$ is finite dimensional.
It can be shown (see e.g. \cite{W1}) that $H_{M+1}$ is spanned linearly by
$(M+1)!$ elements, e.g., those which appear in the expansion of the special operator
$$
\Sigma_{1 \to M+1} =
f_{1 \to M+1} \, f_{1 \to M} \cdots  f_{1 \to 2} \,  f_{1 \to 1} \; ,
$$
where $f_{1 \to n}$ are 1-shuffles defined inductively by $f_{1\to 1}  = 1$,
$f_{1 \to n+1} = 1 + f_{1 \to n} \, \sigma_n$.
Below we assume that $q \neq \exp(2\pi i n/m)$, $n,m \in {\bf Z}$
($q$ is "generic"); for these values of $q$, there exists an isomorphism
between the algebra $H_{M+1}(q)$ and the
group algebra of the symmetric group $S_{M+1}$
(the case $q= \pm 1$ is exceptional, in this case $H_{M+1} =$
group algebra of $S_{M+1}$).

An essential information about a finite dimensional semisimple algebra ${\cal A}$
is contained in the structure of its regular bimodule
which decomposes into direct sums:
${\cal A} = \bigoplus_{\alpha=1}^s  {\cal A} \cdot e_\alpha \; , \;\;\;
{\cal A} =  \bigoplus_{\alpha=1}^s  e_\alpha \cdot {\cal A}$
of left and right submodules (ideals), respectively
(left- and right- Peirce decompositions). Here the elements
$e_\alpha \in {\cal A}$ $(\alpha=1, \dots , s)$ are mutually orthogonal
idempotents: $e_\alpha \, e_\beta = \delta_{\alpha \beta} \, e_\alpha$,
resolving the identity operator:
$1 = \sum_{\alpha=1}^s \, e_\alpha$.
There are two important decompositions of the identity operator
and correspondingly two sets of the idempotents in ${\cal A}$: \\
(1) {\it Primitive idempotents}. An idempotent $e_\alpha$
is {\it primitive} if it can not be further
resolved into a sum of nontrivial mutually orthogonal idempotents. \\
(2) {\it Primitive central idempotents}. An idempotent
$e'_\beta$
is {\it primitive central} if it is primitive in the class
of central idempotents.

%\section{Jucys - Murphy elements}
%\setcounter{equation}0

{} For the A-type Hecke algebra $H_{M+1}(q)$ a set of
elements $\{ y_i \}$ ($i=1, \dots ,M+1$) is defined inductively:
$y_1 =1$, $y_{i+1}= \sigma_i y_i \sigma_i$. These elements
are called {\it Jucys - Murphy elements} and can be written
(using the Hecke condition (\ref{ahecke}) and the braid relation (\ref{braidg}))
in the form
\be
\label{jucmu}
y_i = \sigma_{i-1} \dots \sigma_2 \, \sigma_1^2 \, \sigma_2
\dots \sigma_{i-1} =
(q-q^{-1}) \sum_{k=1}^{i-1}
\, \sigma_{k} \dots \sigma_{i-2} \, \sigma_{i-1} \, \sigma_{i-2}
\dots \sigma_{k} + 1 .
\ee
Sometimes
it is more convenient to use elements $(y_i-1)/(q-q^{-1})$ which,
due to (\ref{jucmu}), have a nontrivial classical limit ($q \to 1$).
The elements $y_i$ pairwise commute. The following statement
explains the importance of the set $\{ y_i \}$.

%\vspace{0.1cm}
\noindent
{\bf Proposition 1.}
{\it The set of Jucys - Murphy elements $\{ y_{i} \}$ ($i=1, \dots ,M+1$)
generates a maximal commutative subalgebra $Y_{M+1}$ in $H_{M+1}$}.

We construct primitive orthogonal
idempotents $e_\alpha \in H_{M+1}$ as functions of
the elements $y_i \in Y_{M+1}$; they are
common {\it eigenidempotents} of $y_i$:
$y_i e_\alpha = e_\alpha y_i  = a^{(\alpha)}_i e_\alpha$ $(i=1,\dots,M+1)$.
We denote (as in \cite{OV}, for symmetric groups)
by ${\rm Spec}(y_1, \dots , y_{M+1})$ the set $\{ \Lambda(e_\alpha) \}$
$(\forall \alpha)$ of strings of eigenvalues:
$\; \Lambda(e_\alpha) =$ $(a^{(\alpha)}_1, \dots , a^{(\alpha)}_{M+1})$.
In view of the following inclusions of the subalgebras $Y_i$ and $H_i(q)$:
$$
\begin{array}{c}
 H_i(q) \subset H_{i+1}(q) \\[-0.1cm]
\!\!\!\! \cup \;\;\;\;\;\;\;\;\; \cup   \\[-0.1cm]
 Y_i  \; \subset \; Y_{i+1}
\end{array}
$$
one can describe the idempotents $\in H_{i+1}$ by considering
the branching of the idempotents of $H_i$ in $H_{i+1}$. It can be shown
that the multiplicity of this branching is equal to one and $y_i$ are
semi-simple for generic $q$.

We need important intertwining operators \cite{Isaev1}
(presented in another form in \cite{Cher5})
\be
\label{impint}
U_{n+1} = \sigma_n y_n - y_n \sigma_n  \;\; (1 \leq n \leq M) \; .
\ee
Elements  $U_{i}$ satisfy relations\footnote{The definition (\ref{impint}) of
intertwining elements is not unique. One can multiply $U_{n+1}$ by
a function $f(y_n, y_{n+1})$: $U_{n+1}  \rightarrow U_{n+1} f(y_n, y_{n+1})$.
Then eqs. (\ref{importt})-(\ref{import}) are valid if
$f(y_n, y_{n+1}) f(y_{n+1}, y_{n}) =1$.}
$U_n \, U_{n+1} \, U_n = U_{n+1} \, U_{n} \, U_{n+1}$ and
\be
\label{importt}
\begin{array}{c}
U_{n+1} y_n = y_{n+1} U_{n+1}  , \;
U_{n+1} y_{n+1} = y_{n} U_{n+1}  ,
\; [U_{n+1}, \, y_k ] = 0 \; (k \neq n,n+1)  ,
\end{array}
\ee
\be
\label{import}
U_{n+1}^2 = (q y_{n} - q^{-1} \, y_{n+1}) \, (q \, y_{n+1} - q^{-1} \, y_n) \; .
\ee
The operators $U_{n+1}$ "permute" elements $y_n$ and $y_{n+1}$
(see (\ref{importt})) which supports a statement
that the center $Z_{M+1}$ of the Hecke algebra
$H_{M+1}$ is generated by symmetric functions in $\{ y_i \}$
$(i = 2, \dots , M+1)$ (to prove this fact it is enough to check relations:
$[\sigma_k, \, y_n + y_{n+1}] = 0 = [\sigma_k, \, y_n  y_{n+1}]$
for all $k < n+1$).

%\vspace{0.1cm}
\noindent
{\bf Proposition 2.} {\it One has
\be
\label{spec}
{\rm Spec}(y_j) \subset \{ q^{2 {\bf Z}_j} \} \;\;\;\; \forall j = 1,2, \dots , M+1 \; ,
\ee
where ${\bf Z}_j$ denotes the set of integers
$\{ 1-j, \dots , -2,-1,0,1,2, \dots ,j-1 \}$}.

%\vspace{0.1cm}
\noindent
{\bf Proof.} We prove (\ref{spec}) by induction.
Obviously, ${\rm Spec} (y_1)$ satisfies (\ref{spec}).
Assume that the spectrum of $y_{j-1}$ satisfies (\ref{spec}) for some $j \geq 2$.
Consider a characteristic equation for $y_{j-1}$ ($j \geq 2$):
$$
f(y_{j-1}) : = \prod_\alpha (y_{j-1} - a^{(\alpha)}_{j-1}) = 0 \;\;\;
(a^{(\alpha)}_{j-1} \in {\rm Spec} (y_{j-1})) \; .
$$
Using properties (\ref{importt})-(\ref{import}) of operators $U_j$, we deduce
\be
\label{spec1}
\begin{array}{c}
0= U_j f(y_{j-1}) U_j = f(y_{j}) U2_j
= f(y_{j}) (q^2 y_{j-1} - y_{j})( y_{j} - q^{-2} y_{j-1}) \; .
\end{array}
\ee
which means that
${\rm Spec} (y_j) \subset
\left( {\rm Spec}(y_{j-1}) \cup q^{\pm 2} \cdot {\rm Spec}(y_{j-1}) \right)$. \hfill $\bullet$

\section{Generalization of the approach of \cite{OV} to the Hecke algebra case}
\setcounter{equation}0

Consider a subalgebra $\hat{H}_2^{(i)}$ in $H_{M+1}$ with generators $y_i$, $y_{i+1}$
and $\sigma_i$ (for fixed $i \leq M$).
We investigate representations of $\hat{H}_2^{(i)}$
with diagonalizable $y_i$ and $y_{i+1}$. Let $e$ be
a common eigenidempotent of $y_i$, $y_{i+1}$:
$y_i e = a_i e$, $y_{i+1} e = a_{i+1} e$. Then the left action of $\hat{H}^{(i)}_2$
closes on elements $v_1 = e$ and $v_2 = \sigma_i e$ and is given by matrices:
\be
\label{abmw1a}
\sigma_i =
\left(
\begin{array}{cc}
0 & 1 \\
1 & q-q^{-1}
\end{array}
\right)  , \;
y_i =
\left(
\begin{array}{cc}
a_i & - (q-q^{-1}) a_{i+1}  \\
0 & a_{i+1}
\end{array}
\right) , \;
y_{i+1} =
\left(
\begin{array}{cc}
a_{i+1} & (q-q^{-1}) a_{i+1}  \\
0 & a_i
\end{array}
\right) ;
\ee
$a_i \neq a_{i+1}$
otherwise $y_i$, $y_{i+1}$ are not diagonalizable.
The matrices $y_i$, $y_{i+1}$ (\ref{abmw1a}) can be simultaneously
diagonalized by a similarity transformation $y \rightarrow V^{-1} y V$, where
$$
V =
\left(
\begin{array}{cc}
1 & \frac{(q-q^{-1}) \, a_{i+1}}{a_i -a_{i+1}} \\
0 & 1
\end{array}
\right)  \; , \;\;\;
V^{-1} =
\left(
\begin{array}{cc}
1 & - \frac{(q-q^{-1}) \, a_{i+1}}{a_i - a_{i+1}} \\
0 & 1
\end{array}
\right)  \; .
$$
As a result we obtain
\be
\label{ah1}
\!\!\! \sigma_i =
\left( \!\!
\begin{array}{cc}
- \frac{(q-q^{-1}) \, a_{i+1}}{a_i -a_{i+1}}  &
1 - \frac{(q-q^{-1})2 \, a_i a_{i+1}}{(a_i -a_{i+1})2} \\[0.3cm]
1  &  \frac{(q-q^{-1}) \, a_i}{a_i -a_{i+1}}
\end{array}
\!\! \right)  , \;
y_i =
\left( \!\!
\begin{array}{cc}
a_i \!\! & \!\! 0  \\
0 \!\! & \!\! a_{i+1}
\end{array}
\!\! \right) , \;
y_{i+1} =
\left( \!\!
\begin{array}{cc}
a_{i+1} \!\! & \!\! 0  \\
0 \!\! & \!\! a_i
\end{array}
\!\! \right) .
\ee
When $a_{i+1} = q^{\pm 2} a_i$,
the 2-dimensional representation (\ref{ah1}) reduces to
a 1-dimensio\-nal representation with $\sigma_i \cdot e = \pm q^{\pm 1} \, e$,
respectively. We summarize the above results as
(cf. Proposition 4.1 \cite{OV}):

\vspace{0.1cm}
\noindent
{\bf Proposition 3.}  {\it Let
$\Lambda = (a_1, \dots , a_i, a_{i+1} , \dots , a_{M+1})
\in {\rm Spec} (y_1, \dots , y_{M+1})$
be a possible spectrum of the set $(y_1, \dots , y_{M+1})$
which corresponds to a primitive idempotent $e_\Lambda \in H_{M+1}$.
Then $a_i = q^{2 m_i}$, where  $m_i \in {\bf Z}_i$ (see Prop. 2) and
(a) $a_i \neq a_{i+1}$ for $i \leq M$;
(b) if $a_{i+1} = q^{\pm 2} a_{i}$ then
$\sigma_i \cdot e_\Lambda = \pm q^{\pm 1} e_{\Lambda}$;
(c) if $a_{i+1} \neq q^{\pm 2} a_{i}$ then
\be
\label{lam'}
\Lambda' = (a_1, \dots , a_{i+1} , a_i, \dots , a_{M+1})
\in {\rm Spec} (y_1, \dots , y_{M+1})
\ee
and the left action of the elements $\sigma_i,y_i,y_{i+1}$
in the linear span of
$v_\Lambda = e_\Lambda$ and
$v_{\Lambda'} = \sigma_i \, e_\Lambda + \frac{(q-q^{-1}) a_{i+1}}{a_i- a_{i+1}} \, e_\Lambda$
is given by (\ref{ah1}).
}

\vspace{0.1cm}
\noindent
{\bf Proposition 4.} {\it Consider the string $\Lambda = (a_1, \dots , a_n)$
of numbers $a_i = q^{2 m_i}$, where  $m_i \in {\bf Z}_i$ (see Prop. 2). Then
$\Lambda = (a_1,a_2, \dots , a_n) \in {\rm Spec}(y_1,y_2, \dots , y_{n})$
iff $\Lambda$ satisfies the following conditions ($z \in {\bf Z}$)
\be
\label{con1}
\!\!\!
\begin{array}{l}
(1) \;\;\;\; a_1 = 1  \; ; \\
(2) \;\;\;\; a_j = q^{2z} \Rightarrow
\{ q^{2(z + 1)},q^{2(z - 1)} \} \cap \{a_1, \dots , a_{j-1} \} \neq \O
%\emptyset
\;\;\; \forall j >1 \; , \;\;  z \neq 0 ; \\
(3) \;\;\;\; a_i = a_j = q^{2z} \; (i < j) \Rightarrow
 \; \{ q^{2(z + 1)},q^{2(z - 1)} \} \subset
\{a_{i+1}, \dots , a_{j-1} \} \; .
\end{array}
\ee
}
\vspace{0.1cm}
\noindent
{\bf Proof.}
The condition (1) is the identity $y_1 =1$.
Conditions (2),(3) can be proven by induction
(see the proof of analogous Theorem 5.1 in \cite{OV}).
To prove the condition (3) we need the
fact that the combinations $(\dots, a_{i-1}, a_i, a_{i+1}, \dots)=$
$(\dots, a, q^{\pm 2} a, a, \dots)$ cannot appear in $\Lambda$:
the braid relation $\sigma_{i} \sigma_{i \pm 1} \sigma_i =
\sigma_{i \pm 1} \sigma_i  \sigma_{i \pm 1}$ is incompatible with
the values $\sigma_i = \pm q^{\pm 1}$, $\sigma_{i+1} = \mp q^{\mp 1}$
(see the condition (b) of Proposition 3).
\hfill $\bullet$

%\vspace{0.5cm}

Consider a Young diagram  with $M+1$ nodes.
We place the numbers $1, \dots , M+1$ into the nodes of the diagram
in such a way that these numbers are arranged along rows and columns in
ascending order in right and down directions. Such diagram
is called a standard Young tableau $[\nu]_{M+1}$.
The standard Young tableau $[\nu]_{M+1}$ defines an ascending set of
standard tableaux:
$[\nu]_{1}  \subset [\nu]_{2}  \subset  \dots  \subset [\nu]_{M+1}$.
In addition we associate a number $q^{2(n-m)}$ (the "content")
to each node of the standard Young tableau, where $(n,m)$ are coordinates of the node.
Example:

\unitlength=6mm
\begin{picture}(25,4)
\put(5,3.4){\vector(1,0){7}}
\put(5,3.4){\vector(0,-1){3.5}}
\put(12,3.5){$n$}
\put(4.2,0){$m$}

\put(5.5,1.3){$
\begin{tabular}{|c|c|c|c|}
\hline
  $\!\!\!\! \!\! ^{1}$  $_1$ & $\!\! ^2$  $_{q^2}$  &
  $\!\! ^4$  $_{q^4}$ & $\!\! ^6$  $_{q^6}$ \\[0.2cm]
\hline
  $\!\! ^3$ $\!\! _{q^{\! -2}}$  &   $\!\!\!\!\! ^5$ $_1$  & $\!\!\! ^8$ $_{q^2}$  &
  \multicolumn{1}{c}{} \\[0.2cm]
\cline{1-3}
  $\!\! ^7$ $\!\! _{q^{\! -4}}$  & \multicolumn{1}{c}{} \\[0.2cm]
\cline{1-1}
\end{tabular}
$}
\end{picture}
\vspace{-1cm}
\be
\label{111}
{}
\ee
\noindent
In general, for the tableau $[\nu]_{M+1}$, the $i$-th node $[\nu]_{i} \backslash [\nu]_{i-1}$
with coordinates $(n,m)$ looks like:
$\begin{tabular}{|c|}
\hline
  $\!\! ^{i}$  $\!\! _{q^{2(n-m)}} \!\!$  \\[0.2cm]
\hline
\end{tabular}$. Thus, to each standard Young tableau $[\nu]_n$ one can associate a
string $(a_1, \dots , a_n)$ with $a_i = q^{2(n-m)}$. E.g.,
a standard Young tableau (\ref{111})
corresponds to a string $(1,q^2,q^{-2},q^4,1,q^6,q^{-4},q^2)$.
This string satisfies conditions of Prop. 3 and therefore
$(1,q^2,q^{-2},q^4,1,q^6,q^{-4},q^2) \in {\rm Spec} (y_1, \dots , y_{8})$.
This relation between contents of $[\nu]_n$ and elements of
${\rm Spec} (y_1, \dots , y_n)$ can be formulated as
(cf. Prop. 5.3 \cite{OV}):

\vspace{0.1cm}
\noindent
{\bf Proposition 5.} {\it There is a bijection between the set $T(n)$ of
the standard Young tableaux with $n$ nodes and the set
${\rm Spec} (y_1, \dots , y_n)$.}

\section{Coloured Young graph and explicit construction of idempotents $e_\alpha$}
\setcounter{equation}0

The above results can be visualized in a different form,
in terms
of a Young graph. By definition, a Young graph is a graph whose
vertices are Young diagrams and edges indicate inclusions of diagrams.
We put the eigenvalues $a_i$ (colours) of the
Jucys-Murphy elements $y_i$  on the edges in such a way
that the string $(a_1, a_2, \dots , a_n)$ along the path from
the top $\emptyset$ of the Young graph to the diagram $\lambda$ with $n$ nodes gives
the content string of the tableau of shape $\lambda$.
For example, the coloured Young graph for $H_4$ is: \\

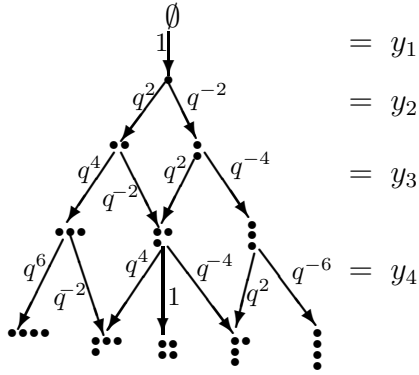
\setlength{\unitlength}{1200sp}%
\begingroup\makeatletter\ifx\SetFigFont\undefined%
\gdef\SetFigFont#1#2#3#4#5{%
  \reset@font\fontsize{#1}{#2pt}%
  \fontfamily{#3}\fontseries{#4}\fontshape{#5}%
  \selectfont}%
\fi\endgroup%
\begin{picture}(6466,7291)(-968,-6894)
\thicklines
{\put(4260,-4336){\vector( 3,-4){1350}}}
{\put(4160,-4336){\vector( 0,-4){1820}}}
{\put(4126,-4411){\vector(-2,-3){1107.692}}}
{\put(2251,-4111){\vector( 1,-3){682.500}}}
{\put(2101,-4186){\vector(-1,-2){930}}
}%
{\put(5026,-2386){\vector( 2,-3){900}}
}%
{\put(4801,-2536){\vector(-1,-2){660}}
}%
{\put(3301,-2386){\vector( 1,-2){770}}}
{\put(3151,-2386){\vector(-2,-3){969.231}}
}%
{\put(4276,-961){\vector( 1,-2){570}}}
{\put(4200,-961){\vector(-3,-4){900}}
}%
{\put(4276,100){\vector( 0,-1){910}}}
{\thinlines
\put(1276,-6136){\circle*{150}}}
{\put(7351,-6811){\circle*{150}}}
{\put(7351,-6586){\circle*{150}}}
\thicklines
{\put(6051,-4486){\vector( -1,-4){380}}}
\put(1901,-5536){\makebox(0,0)[lb]{\smash{\SetFigFont{10}{12.0}{\rmdefault}{\mddefault}
{\updefault}{$q^{\!\! -2}$}%
}}}
\put(3376,-4861){\makebox(0,0)[lb]{\smash{\SetFigFont{10}{12.0}{\rmdefault}{\mddefault}
{\updefault}{$q^4$}%
}}}
\put(4851,-4936){\makebox(0,0)[lb]{\smash{\SetFigFont{10}{12.0}{\rmdefault}{\mddefault}{
\updefault}{$q^{\! -4}$}%
}}}
\put(5551,-2836){\makebox(0,0)[lb]{\smash{\SetFigFont{10}{12.0}{\rmdefault}{\mddefault}
{\updefault}{$q^{-4}$}%
}}}
\put(6826,-5011){\makebox(0,0)[lb]{\smash{\SetFigFont{10}{12.0}{\rmdefault}{\mddefault}
{\updefault}{$q^{-6}$}%
}}}
\put(1276,-4936){\makebox(0,0)[lb]{\smash{\SetFigFont{10}{12.0}{\rmdefault}{\mddefault}
{\updefault}{$q^6$}%
}}}
\put(2401,-2911){\makebox(0,0)[lb]{\smash{\SetFigFont{10}{12.0}{\rmdefault}{\mddefault}
{\updefault}{$q^4$}%
}}}
\put(4256,-5541){\makebox(0,0)[lb]{\smash{\SetFigFont{10}{12.0}{\rmdefault}{\mddefault}
{\updefault}{$1$}%
}}}

\put(8000,-300){\makebox(0,0)[lb]{\smash{\SetFigFont{12}{12.0}{\rmdefault}{\mddefault}
{\updefault}{$= \; y_1$}%
}}}
\put(8000,-1511){\makebox(0,0)[lb]{\smash{\SetFigFont{12}{12.0}{\rmdefault}{\mddefault}
{\updefault}{$= \; y_2$}%
}}}
\put(8000,-3011){\makebox(0,0)[lb]{\smash{\SetFigFont{12}{12.0}{\rmdefault}{\mddefault}
{\updefault}{$= \; y_3$}%
}}}
\put(8000,-5011){\makebox(0,0)[lb]{\smash{\SetFigFont{12}{12.0}{\rmdefault}{\mddefault}
{\updefault}{$= \; y_4$}%
}}}

\put(5876,-5461){\makebox(0,0)[lb]{\smash{\SetFigFont{10}{12.0}{\rmdefault}{\mddefault}
{\updefault}{$q^2$}%
}}}
\put(2906,-3511){\makebox(0,0)[lb]{\smash{\SetFigFont{10}{12.0}{\rmdefault}{\mddefault}
{\updefault}{$q^{\! -2}$}
}}}
\put(4176,-2911){\makebox(0,0)[lb]{\smash{\SetFigFont{10}{12.0}{\rmdefault}{\mddefault}
{\updefault}{$q^2$}%
}}}
\put(4651,-1411){\makebox(0,0)[lb]{\smash{\SetFigFont{10}{12.0}{\rmdefault}{\mddefault}
{\updefault}{$q^{-2}$}%
}}}
\put(3526,-1411){\makebox(0,0)[lb]{\smash{\SetFigFont{10}{12.0}{\rmdefault}{\mddefault}
{\updefault}{$q^2$}%
}}}
\put(4001,-286){\makebox(0,0)[lb]{\smash{\SetFigFont{10}{12.0}{\rmdefault}{\mddefault}
{\updefault}{$1$}%
}}}
{\put(6151,-4411){\vector( 3,-4){1161}}
}%
{\thinlines
\put(7351,-6361){\circle*{150}}
}%
{\put(4200,204){$\emptyset$}}
{\put(4276,-886){\circle*{150}}
}%
{\put(6001,-3886){\circle*{150}}
}%
{\put(6001,-4336){\circle*{150}}
}%
{\put(4276,-4036){\circle*{150}}
}%
{\put(4051,-4036){\circle*{150}}
}%
{\put(4051,-4261){\circle*{150}}
}%
{\put(2476,-4036){\circle*{150}}
}%
{\put(2251,-4036){\circle*{150}}
}%
{\put(2026,-4036){\circle*{150}}
}%
{\put(4876,-2236){\circle*{150}}
}%
{\put(4876,-2461){\circle*{150}}
}%
{\put(3376,-2236){\circle*{150}}
}%
{\put(3151,-2236){\circle*{150}}
}%
{\put(6001,-4111){\circle*{150}}
}%
{\put(1051,-6136){\circle*{150}}
}%
{\put(7351,-6136){\circle*{150}}
}%
{\put(5651,-6736){\circle*{150}}}
{\put(5651,-6511){\circle*{150}}}
{\put(5876,-6286){\circle*{150}}}
{\put(5651,-6286){\circle*{150}}}
{\put(4376,-6586){\circle*{150}}}
{\put(4151,-6586){\circle*{150}}}
{\put(4156,-6361){\circle*{150}}}
{\put(4371,-6361){\circle*{150}}}
{\put(3226,-6286){\circle*{150}}
}%
{\put(3001,-6286){\circle*{150}}
}%
{\put(2776,-6286){\circle*{150}}
}%
{\put(2776,-6511){\circle*{150}}
}%
{\put(1726,-6136){\circle*{150}}
}%
{\put(1501,-6136){\circle*{150}}}

%\put(8401,-7036)
\end{picture}

%\vspace{0.1cm}
\noindent
The path
$\{ \emptyset \stackrel{1}{\rightarrow} \bullet \stackrel{q^2}{\rightarrow} \bullet \, \bullet
\stackrel{q^{-2}}{\rightarrow}
\begin{array}{cc}
\bullet & \!\!\!\! \bullet  \\[-0.25cm]
\bullet
\end{array} \stackrel{1}{\rightarrow}
\begin{array}{cc}
\bullet & \!\!\!\! \bullet  \\[-0.25cm]
\bullet  & \!\!\!\! \bullet
\end{array} \}$
corresponds to the tableau
$[\nu]_4 := \begin{tabular}{|c|c|}
\hline
$\!\!\!^1$ $_1$ $\!\!\!$ & $\!\!\!^2$ $_{q^2} \!\!\!$ \\
\hline
$\!\!\! ^3$ $_{q^{\! - \! 2}} \!\!\!$ & $\!\!\!\! ^4$ $_1 \!\!\!$ \\
\cline{1-2}
\end{tabular}$
with content string $(1,q^2,q^{-2},1)$:
the shape of the tableau is given
by the shape of the last vertex of the path while the labels of nodes
of the tableau shows in which sequence
the points $\bullet$ appear in the vertices along the path. The edge indices of the path
are eigenvalues of the Jucys-Murphy
elements: $(1,q^2,q^{-2},1) \in {\rm Spec}(y_1,y_2,y_3,y_4)$
corresponding to the values of $y_i$ on the primitive idempotent
$e \left( [\nu]_4 \right)$.
Thus, we associate a standard Young tableau with $n$ nodes (related to a string
in ${\rm Spec} (y_1, \dots , y_n)$ and, correspondingly,
to the primitive orthogonal
idempotent of $H_n$) with a path
which starts from the vertex $\emptyset$ and goes down to the vertex
with Young diagram with $n$ nodes (the path with $n$ edges in
the coloured Young graph).
Denote by $X(n)$ the set of all such paths and by ${\rm Str}(n)$ the set of
the strings $\Lambda = (a_1, \dots , a_n)$ of numbers $a_i = q^{2 m_i}$
satisfying conditions (\ref{con1}). We collect the above construction in the
following statement.

\vspace{0.1cm}
\noindent
{\bf Proposition 6.} {\it There is a bijection between the set $T(n)$ of
the standard Young tableaux with $n$ nodes,
the set ${\rm Spec} (y_1, \dots , y_n)$,the set ${\rm Str}(n)$ and the set $X(n)$ of
the paths of length $n$ in the Young graph:
$T(n) \leftrightarrow {\rm Spec} (y_1, \dots , y_n)
\leftrightarrow {\rm Str}(n) \leftrightarrow X(n)$.
}

The dimension of the irreducible representation of $H_n(q)$
(corresponding to the Young diagram $\lambda$ with $n$ nodes)
is equal to the number of standard tableaux $[\nu]_n$ of shape $\lambda$ or, as we saw,
to the number of paths which lead to this Young diagram from the top vertex $\emptyset$.
This number is given by a Frobenius
formula $d_\lambda = n! (h_1! \dots h_k!)^{-1} \prod_{i <j}(h_i - h_j)$, where
$k$ is the number of rows in $\lambda$ and $h_i$ are hook lengths
%(we give the definition of the hook lengths below)
of the nodes
in the first column of $\lambda$ (see, e.g., \cite{OgPya}).

Since the coloured Young graph for $H_{M+1}$ contains
the whole information about the spectrum of $y_k$,
we can deduce the expressions (in terms of the
elements $y_k$) of all orthogonal primitive idempotents
for the Hecke algebra using
the inductive procedure proposed in \cite{OgPya}.
This special set of primitive orthogonal idempotents
has also been described in \cite{Mur}.

Let $\lambda$ be a Young diagram with $n=n_k$ rows:
$\lambda_1 \geq \lambda_2 \geq  \dots \geq \lambda_n$ and
$|\lambda| := \sum_{i=1}^n \lambda_i$ be the number of its nodes.
Consider the case when
$\lambda_1 = \dots = \lambda_{n_1} = \lambda_{(1)} >
\lambda_{n_1 +1} = \lambda_{n_1 + 2} = \dots = \lambda_{n_2} = \lambda_{(2)} >
\dots > \lambda_{n_k - n_{k-1} +1} = \dots = \lambda_{|\lambda|} = \lambda_{(n_k)}$:

\unitlength=4mm
\begin{picture}(25,6.5)
\put(0,3){$\lambda = $}
\put(4,5.5){\line(1,0){5}}
\put(4,4){\line(1,0){5}}
\put(4,4){\line(0,1){1.5}}
\put(9,4){\line(0,1){1.5}}

\put(4,2.5){\line(1,0){3}}
\put(7,2.5){\line(0,1){1.5}}
\put(4,4){\line(1,0){3}}
\put(4,1){\line(0,1){3}}
\put(6,1.5){\line(0,1){1}}
\put(4,1.5){\line(1,0){2}}
\put(4.5,1){$\dots$}

\put(4,0){\line(0,1){0.5}}
\put(5,0){\line(0,1){0.5}}
\put(4,0){\line(1,0){1}}
\put(4,0.5){\line(1,0){1}}

\put(6,6.2){$_{_{\lambda_{_{(1)}}}}$}
\put(3,4.7){$_{n_{_1}}$}
\put(2,3.4){$_{_{n_2-n_1}}$}
\put(1,0.4){$_{_{n_k-n_{k-1}}}$}
\put(9.1,4){$_{n_{_1},} {_{\lambda_{_{(1)}}}}$}
\put(7.1,2.5){$_{n_{_2},} {_{\lambda_{_{(2)}}}}$}
\put(6.2,1.5){$_{n_{_3},} {_{\lambda_{_{(3)}}}}$}
\put(5.1,0){$_{n_{_k},} {_{\lambda_{_{(k)}}}}$}

\end{picture}
\vspace{-1cm}
\be
\label{qdima01}
\ee

\noindent
Here $(n_i,\lambda_{(i)})$ are coordinates of the nodes corresponding
to the corners of the diagram $\lambda$. Consider any
standard Young tableau $[\nu]_{|\lambda|}$ of shape
(\ref{qdima01}). Let $e([\nu]_{|\lambda|}) \in H_{|\lambda|}$
be a primitive idempotent
corresponding to the tableau $[\nu]_{|\lambda|}$.
Taking into account the branching rule implied by the coloured Young graph
for $H_{|\lambda|+1}$ we conclude that the following identity holds
$$
e([\nu]_{|\lambda|}) \prod_{r=1}^{k+1} \left( y_{|\lambda|+1}
- q^{2(\lambda_{(r)} - n_{r-1})} \right) = 0 \; ,
$$
where $\lambda_{(k+1)} = n_0 = 0$.
Thus, for a new tableau $[\nu_j]_{|\lambda|+1}$
which is obtained by adding to the tableau $[\nu]_{|\lambda|}$
of shape (\ref{qdima01}) a new node with coordinates $(n_{j-1} +1, \lambda_{(j)} +1)$
we obtain the following primitive idempotent (after a normalization)
\be
\label{pij}
e([\nu_j]_{_{|\lambda|+1}}) :=
e([\nu]_{_{|\lambda|}}) \prod_{\stackrel{r=1}{_{r \neq j}}}^{k+1} \!
\frac{\left( y_{_{|\lambda|+1}}
- q^{2(\lambda_{(r)} - n_{r-1})} \right)}{\left( q^{2(\lambda_{(j)} - n_{j-1})}
- q^{2(\lambda_{(r)} - n_{r-1})} \right)}  = e([\nu]_{_{|\lambda|}}) \, \Pi_j \,  .
\ee
Using this formula and "initial data" $e \left(
\begin{tabular}{|c|}
\hline
$\!\!\!$  1 $\!\!\!$ \\
\hline
\end{tabular}
\right) = 1$, one can deduce step by step
explicit expressions for all primitive orthogonal idempotents for
Hecke algebras.

\section{q-dimensions for Young diagrams}
\setcounter{equation}0

Consider a linear map $Tr_{d(m+1)}$: $H_{m+1}(q) \to H_m(q)$
from the Hecke algebra $H_{m+1}(q)$ to its subalgebra $H_m(q)$ such that
($\forall X,Y \in H_m(q)$, $Z \in H_{m+1}(q)$)
\be
\label{map}
\begin{array}{c}
Tr_{d(m+1)} ( X ) = z_d \, X  \, , \;\;
Tr_{d(m+1)} ( X \, Z \, Y ) = X \, Tr_{d(m+1)}(Z) \, Y \;\;  \, , \\[0.1cm]
Tr_{d(m+1)} ( \sigma_m^{\pm 1} X \sigma_m^{\mp 1}) =
Tr_{d(m)} (X)  \; , \;\;\;
Tr_{d(m+1)} ( \sigma_m ) = 1 \; , \\[0.1cm]
Tr_{d(m)}  Tr_{d(m+1)} ( \sigma_m Z ) =
Tr_{d(m)}  Tr_{d(m+1)} ( Z \sigma_m) \; ,
\end{array}
\ee
where $z_d$ is a constant which we fix as
$z_d = \frac{1 - q^{-2d}}{q - q^{-1}}$ for later convenience.
Then one can define an Ocneanu's trace
${\cal T}r^{(m+1)}$: $H_{m+1}(q) \to {\bf C}$
as a sequence of maps
${\cal T}r^{(m+1)} := Tr_{d(1)} Tr_{d(2)} \cdots Tr_{d(m+1)}$.

\vspace{0.1cm}
\noindent
{\bf Proposition 7}. {\it Ocneanu's traces of idempotents $e([\nu]_{|\lambda|})$,
$e([\nu']_{|\lambda|})$ corresponding to tableaux
$[\nu]_{|\lambda|}$, $[\nu']_{|\lambda|}$
of the same shape $\lambda$ coincide. Thus,
$$
{\rm qdim}(\lambda) := {\cal T}r^{(|\lambda|)} e([\nu]_{|\lambda|}) =
{\cal T}r^{(|\lambda|)} e([\nu']_{|\lambda|})  \;
$$
depends on the diagram $\lambda$ only.
}

Using (\ref{map}) we deduce an identity (see Appendix)
\be
\label{gener1a'}
1 + (q- q^{-1}) Tr_{_{d(|\lambda|+1)}} \left(  \frac{y_{|\lambda|+1} \,
\tau}{1- y_{|\lambda|+1} \, \tau} \right)
= \frac{(1- \tau \, q^{-2d})}{(1-\tau)} \,
\prod_{k=1}^{|\lambda|} \, \frac{(1-\tau \, y_{k})^2}{(1-q^2 \tau y_{k}) (1-q^{-2}\tau y_{k})} \; ,
\ee
where $\tau$ is a parameter. To calculate "qdim" for the diagram (\ref{qdima01})
we need to find the value of the element
(\ref{gener1a'}) on the idempotent $e([\nu]_{|\lambda|})$,
where $[\nu]_{_{|\lambda|}}$ is any Young tableau of shape (\ref{qdima01}).
We take the "row-standard" tableau $[\nu]_{_{|\lambda|}}$
corresponding to the eigenvalues of $y_k$ arranged along the rows from left to right
and from top to bottom:
$$
\begin{array}{l}
y_1 = 1, \; y_2 = q^2, \; y_3 = q^4, \; \dots , \; y_{\lambda_1-1} = q^{2(\lambda_1 -2)},
\; y_{\lambda_1} = q^{2(\lambda_1 -1)}, \\
y_{\lambda_1 +1} = q^{-2}, \; y_{\lambda_1 +1} = 1, \;
\dots , \; y_{\lambda_1 +\lambda_2} = q^{2(\lambda_2 -2)}, \\
\dots\dots\dots\dots\dots\dots , \\
y_{|\lambda|- \lambda_n +1} = q^{-2(n-1)}, \; \dots , \; y_{|\lambda|} = q^{2(\lambda_n -n)} \; .
\end{array}
$$
The result is $(n_k =n, n_0 := 0)$
\be
\label{qdim11}
\!\! Tr_{\!\!_{d(|\lambda|+1)}} \! \left( \! \sum_j  P_{j}
\frac{(q-q^{-1}) \, \mu_j \, \tau}{
1- \mu_j \, \tau} \!\! \right) =
e([\nu]_{|\lambda|}) \! \left( \! \frac{1- \tau \, q^{-2d}}{1-\tau q^{-2n}} \,
\prod_{j=1}^k \frac{1-\tau \, q^{2(\lambda_{(j)}-n_j)}}{
1-\tau \,  \mu_j} - 1 \!\! \right) ,
\ee
where we have inserted into the l.h.s.  the spectral decomposition of
the idempotent $e([\nu]_{_{|\lambda|}})$ (see (\ref{pij})):
$$
e([\nu]_{_{|\lambda|}}) = e([\nu]_{_{|\lambda|}}) \sum_j  \Pi_j = \sum_j P_{j} \; , \;\;\;
P_{j} \, y_{_{|\lambda|+1}} =
P_{j} \, q^{2(\lambda_{(j)} -n_{j-1})}
 = P_{j} \, \mu_j \; .
$$
The operator $P_j$ projects $y_{|\lambda|+1}$ on its eigenvalue
$\mu_j := q^{2(\lambda_{(j)} -n_{j-1})}$
which appeared in the denominator of the r.h.s. of (\ref{qdim11}).
Comparing both sides of eq. (\ref{qdim11}) we deduce
$$
Tr_{\!_{d(|\lambda|+1)}} \left( P_j \right)
=e([\nu]_{_{|\lambda|}}) \,
 \lim_{\;\; \tau \to 1/\mu_j}
\frac{(1- \mu_j \, \tau)}{(q-q^{-1})}
 \left( \frac{1- \tau \, q^{-2d}}{1-\tau q^{-2n}} \,
\prod_{r=1}^k \frac{1-\tau \, q^{2(\lambda_{(r)}-n_r)}}{
1-\tau \, \mu_r } \right)
$$
\be
\label{qdim13}
=e([\nu]_{|\lambda|}) \cdot q^{-d} \, [q^{(\lambda_{(j)}-n_{j-1}+d)}]_q \,
 \frac{ \prod_{n,m \in \lambda}  [h_{n,m}]_q}{
 \prod_{n,m \in \lambda^{(j)}}  [h_{n,m}]_q} \; ,
\ee
where $h_{n,m}$ are hook lengths
of nodes $(n,m)$ of the diagrams $\lambda$ or $\lambda^{(j)}$
($\lambda^{(j)}$ is a diagram obtained by adding to the diagram $\lambda$
a new node with coordinates $(n_{j-1} +1, \lambda_{(j)} +1)$).
Applying the Ocneanu's
trace ${\cal T}r^{(|\lambda|)}$
to eq. (\ref{qdim13}) we find a recurrent relation:
$$
{\rm qdim}(\lambda^{(j)}) = {\rm qdim}(\lambda) \,
q^{-d} \, [\lambda_{(j)}-n_{j-1}+d]_q \,
 \frac{ \prod_{n,m \in \lambda}  [h_{n,m}]_q}{
 \prod_{n,m \in \lambda_j}  [h_{n,m}]_q} \; ,
$$
which is solved by
$$
{\rm qdim}(\lambda) = q^{-d |\lambda |} \,
 \prod_{n,m \in \lambda} \frac{[d + m-n]_q}{[h_{n,m}]_q} \; .
 $$
Up to a normalization factor this formula has firstly been obtained in \cite{W1}.

{}For $R$-matrix representations of $H_{M+1}(q)$
(about $R$-matrix representations of the Hecke algebra see Refs.
\cite{Jimb1}, \cite{10}) which corresponds to
the quantum supergroup $GL_q(N|M)$, the parameter $d$ equals
$N-M$. This justifies our
choice of the parametrization of $z_d$ in the first eq. of (\ref{map}).

Proposition 7 can be generalized. Let $T$ be a quantum matrix satisfying
\be
\lb{qgr}
\hat{R}_{12} \, T_1 \, T_2 = \hat{R}_{12} \, T_1 \, T_2
\ee
in the notations of \cite{10},
where $\hat{R}_{12} = \rho(\sigma_1)$ is the $R$-matrix representation
of the Hecke algebra.

\vspace{0.3cm} \noindent {\bf Proposition 8.} {\it The quantum traces 
(for the definition of the quantum trace see e.g. 
\cite{10}, \cite{11}, \cite{12}) of the matrices 
$[T_1 \cdots T_{|\lambda |} \, \rho (e([\nu]_{|\lambda|}))]$ and 
$[T_1 \cdots T_{|\lambda |} \, \rho (e([\nu']_{|\lambda|}))]$ 
$$ 
\chi_\lambda(T) := Tr_{R(1 \to |\lambda |)} \left( T_1 \cdots T_{|\lambda |} \, 
\rho (e([\nu]_{|\lambda|})) \right) = Tr_{R(1 \to |\lambda |)} \left( T_1 \cdots T_{|\lambda |} \, 
\rho (e([\nu']_{|\lambda|})) \right) \; ,
 $$ 
corresponding to tableaux $[\nu]_{|\lambda|}$ and $[\nu']_{|\lambda|}$ of the same shape $\lambda$, 
coincide. Thus, $\chi_\lambda(T)$ depends only on the diagram $\lambda$. }

Consider the $GL_q(N)$ quantum group (\ref{qgr}) with a standard $GL_q(N)$ 
Drinfeld-Jimbo $R$-matrix $\hat{R}_{12}$ \cite{10}. It is known \cite{Jimb1}, \cite{10} 
that the standard $GL_q(N)$ matrix $\hat{R}_{12}$ defines the representation of the 
Hecke algebra. We note that the $GL_q(N)$ quantum matrix $T$ can be realized by arbitrary 
numerical diagonal $(N \times N)$ matrix $X$. Then $\chi_\lambda(X)$ is a numerical 
function of the deformation parameter $q$ and the entries of $X$. In the classical 
limit $q \to 1$ the operator $\rho (e([\nu]_{|\lambda|}))$ tends to the Young 
projector and the function $\chi_\lambda(X)$ coincides with a character of the 
element $X$ $(X \in GL(N))$ in the representation corresponding to the diagram $\lambda$.

\section{Appendix}
\setcounter{equation}0

Taking into account the definition of the generators $y_m$ we have equations
\be
\lb{qdim3}
\frac{1}{(t-y_{m+1})} \sigma^{-1}_m = \sigma^{-1}_m \frac{1}{(t-y_m)} +
\frac{(q-q^{-1}) y_m}{(t-y_{m+1})} \frac{1}{(t-y_m)} \ee \be \lb{qdim4}
 \frac{1}{(t-y_{m+1})} \sigma_m = \sigma^{-1}_m  \frac{1}{(t-y_m)} +
\frac{ (q-q^{-1}) t}{(t-y_m)} \frac{1}{(t-y_{m+1})} \; .
\ee
Eqs. (\ref{qdim3}),  (\ref{qdim4}) and the definition of the map (\ref{map}) give
a recurrent relation
\be
\lb{qdim7}
\frac{(t-q^2 y_{m}) (t-q^{-2} y_{m})}{(t-y_{m})^2}  Z_{m+1}  = Z_m +
\frac{(q-q^{-1}) y_m}{(t-y_{m})^2}
\left[ 1 - (q-q^{-1}) z_d \right] \; ,
\ee
where the parameter $z_d$ is introduced in (\ref{map}) and
$$
Z_m : = Tr_{d (m)}\left( \frac{1}{(t-y_m)} \right)  \; .
$$
Eq. (\ref{qdim7}) is simplified by the substitution
$Z_m = \tilde{Z}_m -  \left[ 1 - (q-q^{-1}) z_d \right]/((q-q^{-1}) t)$ and we have
$$
\frac{(t-q^2 y_{m}) (t-q^{-2} y_{m})}{(t-y_{m})^2}  \tilde{Z}_{m+1}  = \tilde{Z}_m  \; .
$$
This equation can be easily solved and finally we obtain the expression
$$
Z_{m+1} = \frac{1}{(q-q^{-1}) t}  \left( 1 + \frac{(q-q^{-1}) \, z_d}{(t-1)} \right)
\prod_{k=1}^m \, \frac{(t-y_{k})^2}{(t-q^2 y_{k}) (t-q^{-2} y_{k})} 
$$
$$
- \frac{1}{(q-q^{-1}) t} \left[ 1 - (q-q^{-1}) z_d \right] \; ,
$$
which is equivalent to (\ref{gener1a'}) for $t = 1/\tau$.

\end{document}